\documentclass[11pt]{article}

\usepackage{psfig}

\def\proof{\bf \medbreak \noindent Proof. \rm} 
\def\eoproof{{\unskip\nobreak\hfil\penalty50
	\hskip2em\hbox{}\nobreak\hfil\vrule height4pt width5.5pt depth2pt
	\parfillskip=0pt\finalhyphendemerits=0\medbreak}}

\newcommand{\be}{\begin{equation}}
\newcommand{\ee}{\end{equation}}

\newcommand{\bea}{\begin{eqnarray}}
\newcommand{\eea}{\end{eqnarray}}

\newcommand{\beas}{\begin{eqnarray*}}
\newcommand{\eeas}{\end{eqnarray*}}

\def\CC{\hbox{\rlap{$\,\,
  $\hbox{\vrule height6.2pt width.35pt depth-0.1pt}}$\rm C$}}

\begin{document}
 \title{The Askey Scheme for Hypergeometric Orthogonal 
Polynomials Viewed from Asymptotic Analysis}

\author{Nico M. Temme\\
  CWI\\
  P.O. Box 94079\\
  1090 GB Amsterdam\\
  The Netherlands\\
e-mail: {\tt  nicot@cwi.nl}
 \and
  Jos\'e L. L\'opez\\
Dpto. Matematica e Informatica\\
Universidad Publica de Navarra\\
Campus de Arrosadia s/n\\
31006-Pamplona\\ 
Spain\\
e-mail: {\tt jl.lopez@unavarra.es}}

\maketitle

 \begin{abstract}
Many limits are known for hypergeometric orthogonal polynomials
that occur in the Askey scheme. We show how asymptotic representations
can be derived by using the generating functions of the polynomials. For 
example, we discuss the asymptotic representation of the
Meixner-Pollaczek, Jacobi, Meixner, and
Krawtchouk polynomials in terms of Laguerre polynomials. 
 \end{abstract}




\section{Introduction}\label{intro}

It is well known that the Hermite polynomials
play a crucial role in
certain limits of the classical orthogonal polynomials. For example,  
the ultraspherical (Gegenbauer) polynomials
$C_n^{\gamma}(x)$, which are defined by the generating function
\be
(1-2xw+w^2)^{-\gamma}=\sum_{n=0}^\infty 
C_n^{\gamma}(x) w^n,\quad -1\le x\le 1,\quad |w|<1,\label{i1}
\ee
have the well-known limit
\be
\lim_{\gamma\to\infty}\gamma^{-n/2}
C_n^\gamma(x/\sqrt{{\gamma}})=\frac{1}{n!}H_n(x).\label{i2}
\ee

For the Laguerre polynomials, which are defined by the generating
function
\be
(1-w)^{-\alpha-1}e^{-wx/(1-w)}=\sum_{n=0}^\infty
L_n^\alpha(x)\,w^n,\quad|w|<1,\label{i3}
\ee
$\alpha,x\in\CC$, a similar results reads
\be
\lim_{\alpha\to\infty}\alpha^{-n/2}
L_n^\alpha\left(x\sqrt{{\alpha}}+\alpha\right)=
\frac{(-1)^n\,2^{-n/2}}{n!}
\,H_n\left(x/\sqrt{{2}}\right).\label{i4}
\ee

These limits give insight in the location of the zeros for large
values of the limit parameter, and the asymptotic relation with the Hermite
polynomials if the parameters $\gamma$ and $\alpha$ become large and $x$ is
properly scaled.

Many methods are available to prove these and other limits.
In this paper we concentrate on asymptotic relations between the 
polynomials,
from which the limits may follow as special cases.

In \cite{koeksw} many relations are given for
hypergeometric orthogonal polynomials and their $q-$analogues, 
including limit 
relations between many polynomials. In Figure 1 we show examples for 
which limit relations between neighboring polynomials are available, 
but many other limit relations are mentioned in \cite{elblaf}, 
\cite{grza}, \cite{koeksw} and \cite{rzag}. 

In \cite{loptem1}, \cite{loptem2} and \cite{loptem3} 
we have given several asymptotic
relations between polynomials and Hermite polynomials. 
In these first papers we
considered 
Gegenbauer  polynomials, 
Laguerre  polynomials, 
Jacobi  polynomials, 
Tricomi-Carlitz  polynomials,
generalized Bernoulli polynomials,
generalized Euler  polynomials,
generalized Bessel polynomials
and
Buchholz polynomials.

The method for all these cases is the same and we observe that the 
method also works for polynomials outside the class of hypergeometric
polynomials, such as Bernoulli and Euler  polynomials. 

Our method is different from the one described in 
\cite{elblaf} and
\cite{grza}, where also more terms in the limit relation are
constructed in order to obtain more insight in the limiting process.
In these papers  expansions of the form
$$P_n(x;\lambda)=\sum_{k=0}^\infty R_k(x;n)\lambda^{-k}$$
are considered, which generalizes the limit relation
$$\lim_{\lambda\to\infty}P_n(x;\lambda)=R_0(x;n),$$ 
and which gives deeper information on
the limiting process.  In \cite{grza} a method 
for the recursive computation of the coefficients $R_k(x;n)$ is designed.

In \cite{rzag} similar methods are used,
now in particular for limits between classical discrete (Charlier, Meixner, 
Krawtchouk, Hahn) to classical continuous (Jacobi, Laguerre, Hermite)
orthogonal polynomials.

\vspace*{0.3cm}
\centerline{\protect\hbox{\psfig{file=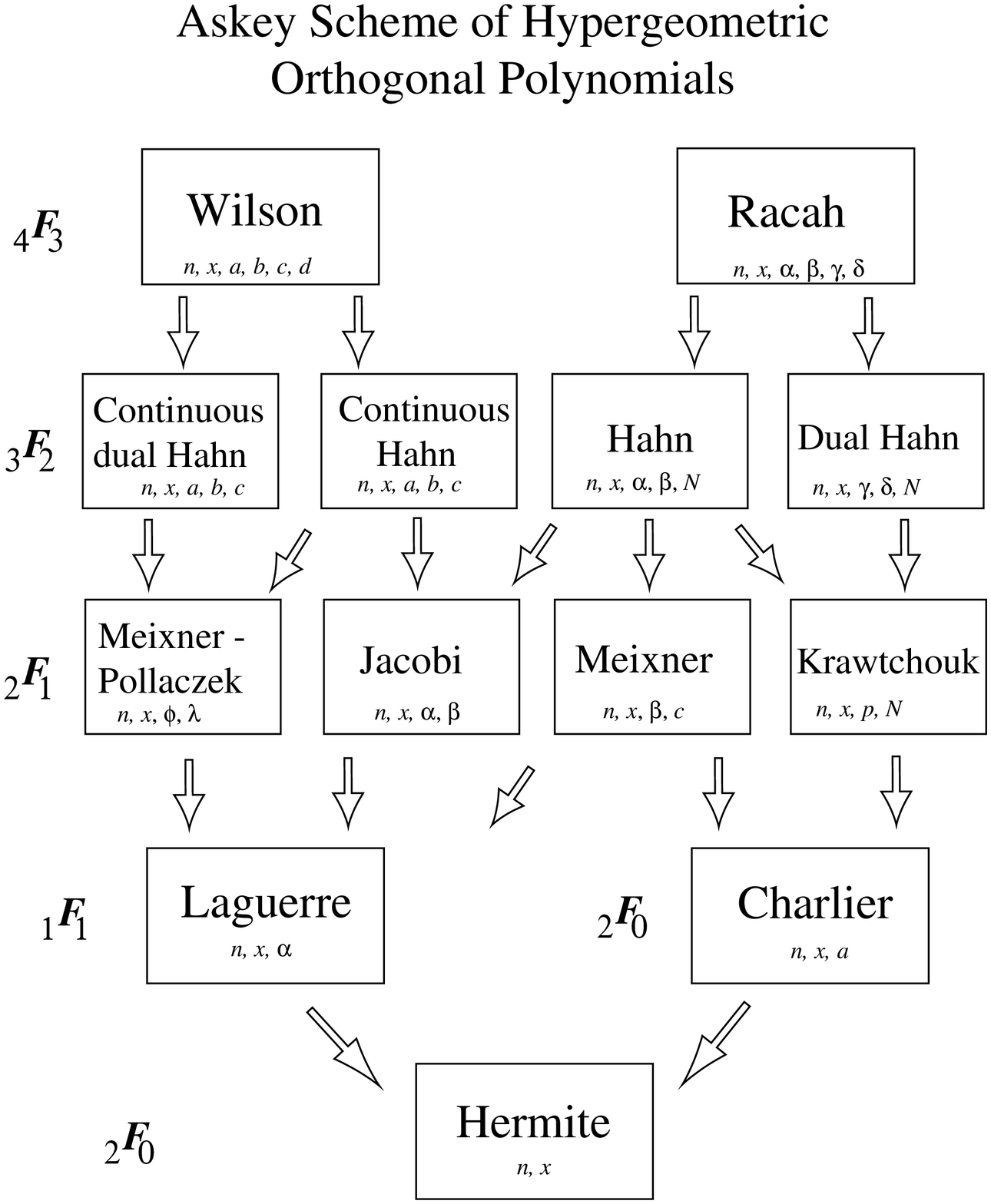,width=12cm}}}
\noindent
{\bf Figure 1.}\quad
The Askey scheme for hypergeometric orthogonal polynomials, with indicated
limit relations between the polynomials.

In current research we investigate if other limits 
in the Askey scheme can be
replaced by asymptotic results. 
Until now we verified all limits from the third
level to the fourth (Laguerre and Charlier) and the fifth level (Hermite). 
Several limits are new, and all results have full asymptotic expansions.

\section{Asymptotic representations}\label{sectwo}

Starting point in our method is a generating series
\be F(x,w)=\sum_{n=0}^\infty p_n(x)\,w^n, \label{t1}
\ee
$F$ is a given function, which is analytic with respect to $w$ at $w=0$, and
$p_n$ is independent of $w$. 

The relation (\ref{t1}) gives for $p_n$ the
Cauchy-type integral
$$p_n(x)=\frac1{2\pi i} \int_{{\cal C}} F(x,w)\,\frac{dw}{w^{n+1}},$$
where ${{\cal C}}$ is a circle around the origin inside the domain where $F$ is
analytic (as a function of $w$).

We write
$$F(x,w)=e^{Aw-Bw^2}\,f(x,w),$$
where $A$ and $B$ do not depend on $w$.
This gives
\be
p_n(x)=\frac1{2\pi i} \int_{{\cal C}} e^{Aw-Bw^2}
\,f(x,w)\,\frac{dw}{w^{n+1}}.\label{t2}
\ee
Because $f$ is also analytic (as a function of $w$), we can expand
\be
f(x,w)=e^{-Aw+Bw^2}F(x,w)=\sum_{k=0}^\infty c_k w^k, \label{t3}
\ee
that is,
$$f(x,w)=1+[p_1(x)-A]w+\left[p_2(x)-Ap_1(x)+
B+\frac12A^2\right]w^2+\ldots$$
if we assume that $p_0(x)=1$ (which implies $c_0=1$).

We substitute (\ref{t3}) in
(\ref{t2}). The Hermite polynomials have the generating function
$$
e^{2xw-w^2}=\sum_{n=0}^\infty\frac{H_n(x)}{n!}w^n,\quad x,w\in\CC,
$$
which gives the Cauchy-type integral
\be
H_n(x)=\frac{n!}{2\pi i}\, \int_{{\cal C}} 
e^{2xz-z^2} z^{-n-1}\,dz,\label{t4}
\ee
where ${{\cal C}}$ is a circle around the origin and the integration is in
positive direction. The result is the finite expansion
\be
p_n(x)=z^n \, \sum_{k=0}^n \frac{c_k}{z^k}\,\frac{H_{n-k}(\xi)}{(n-k)!},
\quad z=\sqrt{{B}}, \quad \xi=\frac{A}{2\sqrt{{B}}},\label{t5}
\ee
because terms with $k>n$ do not contribute in the integral
in (\ref{t2}).

In order to obtain an asymptotic property of (\ref{t5}) we take $A$ and
$B$ such that $c_1=c_2=0$. This happens if we take
$$A=p_1(x),\quad B=\frac12p_1^2(x)-p_2(x).$$

As we will show, the asymptotic property follows from the behavior of the
coefficients $c_k$ if we take a parameter of the polynomial $p_k(x)$ large.
We use the following lemma, and explain what happens by considering a few
examples.

\medskip\noindent
{\bf Lemma 2.1\quad}
\it
Let $\phi(w)$ be analytic at $w=0$, with Maclaurin expansion of the
form 
$$\phi(w)=\mu w^n(a_0+a_1w+a_2w^2+\ldots),$$
where $n$ is a positive integer and $a_k$ are 
complex numbers that do not depend on the complex number $\mu$, $a_0\ne0$.
Let $c_k$ denote the coefficients of the power series of $f(w)=e^{\phi(w)}$,
that is, 
$$f(w)=e^{\phi(w)}=\sum_{k=0}^\infty c_kw^k.$$
Then  $c_0=1, c_k=0, k=1,2\ldots, n-1$ and
$$c_k={{\cal O}}\left(\vert\mu\vert^
{\lfloor k/n\rfloor}\right),\quad \mu\to\infty.$$
\rm
\proof
The proof follows from expanding 
\beas
\sum_{k=0}^\infty c_kw^k&=&e^{\phi(w)}=
\sum_{k=0}^\infty \frac{[\phi(w)]^k}{k!}\\
&=&\sum_{k=0}^\infty \frac{\mu^k w^{kn}}{k!}(a_0+a_1w+a_2w^2+\ldots)^k,
\eeas 
and comparing equal powers of $w$.
\eoproof
            
\subsection{Ultraspherical polynomials}

The generating function is given in (\ref{i1}), and we obtain
$$A=C_1^\gamma(x)=2x\gamma,\quad 
B=\frac12\left[C_1^\gamma(x)\right]^2-C_2^\gamma(x)=\gamma(1-2x^2).$$
The expansion reads
\be
C_n^{\gamma}(x)=z^n\,\sum_{k=0}^n\,\frac
{c_k}{z^k}\,\frac{H_{n-k}(\xi)}{(n-k)!},\label{t6}
\ee
where $ z=\sqrt{{\gamma(1-2x^2)}}, \xi=x\gamma/z$. We have
$$c_0=1, \quad c_1=c_2=0, \quad c_3= \frac23\gamma x(4x^2-3).$$
Higher coefficients follow from a recursion relation.

The function $f(x,w)$ of (\ref{t3}) has the form 
$$f(x,w)=e^{\phi(x,w)}, \quad 
\phi(x,w)=\gamma w^3(a_0+a_1w+a_2w^2+\ldots).
$$
By using Lemma 2.1 and $\xi={{\cal O}}
(\sqrt{{\gamma}})$ we conclude that the
sequence
$\{\phi_k\}$ with $\phi_k=c_k/z^kH_{n-k}(\xi)$ has the following
asymptotic property:
$$\phi_{k}={{\cal O}}\left(\gamma^{n/2+\lfloor k/3\rfloor-k}\right),
\quad k=0,1,2,\ldots.$$
This explains the asymptotic nature
of the representation in (\ref{t6}) for large values of $\gamma$,
with $x$ and $n$ fixed. 

To verify the limit given in (\ref{i2}), 
we first write $x$ in terms of $\xi$: 
$x={\xi}/{\sqrt{{\gamma+2\xi^2}}}$.
With this value of $x$ we can verify that  
$c_k/z^k=o(1),  \gamma\to\infty$, 
and in fact we have the limit
$$
\lim_{\gamma\to\infty}\frac{\gamma^n}{(\gamma+2x^2)^{n/2}}\,
C_n^\gamma\left(\frac{x}{\sqrt{{\gamma+2x^2}}}\right)
=\frac{1}{n!}H_n(x).$$

\subsection{Laguerre polynomials}

We take as generating function (see (\ref{i3}))
$$F(x,w)=(1+w)^{-\alpha-1}e^{wx/(1+w)}=
\sum_{n=0}^\infty (-1)^nL_n^{\alpha}(x)\,w^n.
$$
We have $A=x-\alpha-1,\quad B=x-\frac12(\alpha+1)$, and we obtain
\be
L_n^{\alpha}(x)=(-1)^n\,z^n\,\sum_{k=0}^n\,\frac {c_k}{z^k}\,
\frac{H_{n-k}(\xi)}{(n-k)!},\label{l1}
\ee
where $ z=\sqrt{{x-(\alpha+1)/2}}, \xi=({x-\alpha-1})/({2z})$.
The first coefficients are
$$c_0=1, \quad c_1=c_2=0, \quad c_3= \frac13(3x-\alpha-1).$$
Higher coefficients follow from a recursion relation.
The representation in (\ref{l1}) has an asymptotic character
for large values of $|\alpha|+|x|$.
It is not difficult to verify that the limit given 
in (\ref{i4}) follows from (\ref{l1}).

\section{Expansions in terms Laguerre polynomials}\label{secthree}

We give examples on how to use Laguerre polynomials for approximating 
other polynomials. The method for the Hermite polynomials demonstrated in the
previous section can be used in a similar way. 

\medskip\noindent
{\bf Lemma 3.1\quad}
\it
Let the polynomials $p_n(x)$ be defined by the generating function
$$ F(x,w)=\sum_{n=0}^\infty p_n(x)\,w^n, $$
where $F(x,w)$ is analytic in $w=0$ and  $F(x,0)=1$. 
Let 
$$f(x,w)=e^{-Aw/(Bw-1)}(1-Bw)^{C+1}\,F(x,w),$$ 
and let the coefficients $c_k(x)$ be
defined by the expansion
\be
f(x,w)  = \sum_{k=0}^\infty c_k(x) w^k,\quad c_0=1,\label{f1}
\ee
where $A, B$ and $C$ do not depend on $w$. Then $p_n(x)$ can be represented as
the finite sum 
\be
p_n(x)=B^{n/2} \, \sum_{k=0}^n
{c_k(x)\over B^{k/2}} \,L_{n-k}^{(C)}(\xi),
\quad \xi={A\over B}, \label{f2}
\ee
where $L_n^{\alpha}(x)$ are the Laguerre polynomials. 

\rm
\proof
The polynomials $p_n(x)$ can be written as
$$p_n(x)=\frac1{2\pi i} \int_{{\cal C}} 
e^{Aw/(Bw-1)}(1-Bw)^{-C-1}\, f(x,w) \frac{dw}{w^{n+1}},$$
where ${{\cal C}}$ is a circle around 
the origin in the domain where $F(x,w)$ is
analytic (as a function of $w$). By substituting the expansion of $f(x,w)$ and
using the generating function (\ref{i3}) of the Laguerre polynomials the proof
follows.
\eoproof

This time,  $A, B$ and $C$ can be chosen 
such that $c_1=0, \ c_2=0, \ c_3=0$.
These coefficients are given by
\beas
c_1  &=& p_1 -B C-B+A,\\
c_2  &=& p_2-p_1 BC-p_1 B+p_1  A -A B C+\frac12(B^2 C^2+B^2 C+A^2),\\
c_3  &=&  
p_3-p_2  B C-p_2  B+p_2  A-p_1  A B C-\frac16 (B^3 C^3+ B^3 C+ A^3)+\\
\quad&{\ }&\quad\frac12 (p_1  B^2 C^2+ p_1  
B^2 C+ p_1  A^2 + A B^2 C^2- A B^2 C- B A^2
C+ A^2 B). 
\eeas
We see that the equations $c_1=0, \ c_2=0, \ c_3=0$
for solving for $A, B$ and $C$ are nonlinear. 
However, solving $c_1=0, c_2=0$ for
$A$ and $C$ gives
$$A=B(C+1)-p_1,\quad C=\frac{p_1^2-2p_2+2p_1B-B^2}{B^2},$$
and with these values $c_3$ becomes
$$ c_3 = p_3-p_2 p_1 + \frac13 (p_1 B^2  +  p_1^3  
+ 2 p_1^2  B - 4 p_2 B) ,$$
and $c_3=0$  is a quadratic equation for $B$.

As follows from the above representation of $C$, 
this quantity will depend on $x$. This
gives an expansion for $p_n(x)$ in terms of Laguerre 
polynomials $L_k^C(\xi)$ with 
the order depending on $x$. When studying properties 
of $p_n(x)$ (for example
investigating the zeros) this may not be very desirable. 
In that case we can always take
$C=\alpha$ (not depending on $x$), and concentrate 
on two equations $c_1=0, \ c_2=0$ for
solving $A$ and $B$. This gives
\be
A=\sqrt{p_1^2-(\alpha+1)(2p_2-p_1^2)},\quad 
B=\frac{p_1+A}{\alpha+1}.\label{f4}
\ee
The order $\alpha$ may be chosen conveniently, without requiring $c_3=0$.

{}For large values of certain parameters in $p_n(x)$ expansion (\ref{f2}) 
may have an asymptotic property when taking $c_1= c_2= c_3=0$, but also
when only $c_1=0$ or $c_1=c_2=0$. In the following section we give 
four examples, for one level of the Askey scheme, namely for the 
Meixner-Pollaczek, Jacobi, Meixner, and Krawtchouk polynomials. 

\def\asp{r}
\def\pol{P_n^{(\lambda)}(x;\phi)}
\def\WW{\sqrt{x^2+\lambda^2}}
\section{Expanding Meixner-Pollaczek into 
Laguerre po\-lynomials}\label{sectfour}

{}For the Meixner-Pollaczek polynomials we have the generating function:
\be F(x,w)=\left(1-e^{i\phi}w\right)^{-\lambda+ix}
\left(1-e^{-i\phi}w\right)^{-\lambda-ix}=
\sum_{n=0}^\infty \pol w^n.\label{f5}
\ee
The expansion for the Meixner-Pollaczek polynomials reads
\be
\pol=
\sum_{k=0}^n B^{n-k}c_k
L_{n-k}^{(C)}(\xi),\quad \xi=A/B, \label{f6}
\ee
where the coefficients $c_k$ follow from 
(\ref{f1}) with $F(x,w)$ given in (\ref{f5}).

We write $x+i\lambda=r e^{i\theta},  \theta\in[0,\pi], r\ge0$,
and consider $r\to\infty$; the asymptotic results hold uniformly
with respect to $\theta$.

\subsection{One free parameter}

{}First we consider a simple case by taking $B=1$ and
$C=\alpha$, and solve $c_1=0$ for $A$. This gives
$$
A=\alpha+1-2\lambda\cos\phi-2x\sin\phi.
$$
The first coefficients
$c_k$ are given by
$$
c_0=1,\quad c_1=0,\quad
c_2=x\sin2\phi+\lambda\cos2\phi-2(x\sin\phi+
\lambda\cos\phi)+{1\over2}\alpha,
$$
and the remaining ones can be obtained from the recursion
\bea
(k+1)c_{k+1}&=&2(1+\cos\phi)kc_k+\nonumber\\
&{}&[\alpha+1-2\lambda+4(\cos\phi-1)(\lambda\cos\phi+x\sin\phi)+\nonumber\\
&{}&2(1-k)(1+2\cos\phi)]c_{k-1}+\\
&{}&[4\lambda+2(1+\cos\phi)(k-2)-2(\alpha+1)\cos\phi]c_{k-2}+\nonumber\\
&{}&(\alpha+4-k-2\lambda)c_{k-3}.\nonumber\label{f7}
\eea
The asymptotic property follows from the fact that, 
as in Lemma 2.1, the function
$f(x,w)$ can be written as $f(x,w)=\exp[\psi(w)]$, where
$\psi(w)=rw^2(a_0+a_1w+\ldots)$. Hence, the coefficients 
$c_k$ have the asymptotic
behaviour $c_k={{\cal O}}(r^{\lfloor k/2\rfloor})$, 
as $r\to\infty$. The first term approximation
can be written as
$$
\pol=
L_{n}^{(\alpha)}(\xi)+{{\cal O}}\left(\asp^{n-1}\right),
\quad \xi=A,\quad r\to\infty.
$$

In this case a limit for large values of $r$ 
(or $\lambda$ or $x$) cannot obtained from
the above representations. We can obtain a limit by putting
$\lambda=(\alpha+1)/2$. 
Then, as follows from the recursion relation (\ref{f7}), we have
$c_k={\cal O}(\phi^2)$ as $\phi\to 0$, and we obtain the 
limit of the Askey scheme 
$$
\lim_{\phi\to 0}P_n^{(\alpha+1)/2}
\left[(\alpha+1)(1-\cos\phi)-\xi)/(2\sin\phi);\phi\right]=
L_{n}^{(\alpha)}(\xi).
$$
This includes the limit of the Askey scheme (cf. \cite{koeksw})
$$
\lim_{\phi\to 0}P_n^{(\alpha+1)/2}
(-\xi/(2\phi);\phi)=L_{n}^{(\alpha)}(\xi).
$$

\subsection{Two free parameters}

Next we solve $c_1=0, c_2=0$ for $A$ and $B$, 
with  $C=\alpha$. This gives (cf. (\ref{f4}))
\beas
A&=&\sqrt{4(\lambda\cos\phi+x\sin\phi)^2-
2(\alpha+1)(\lambda\cos2\phi+x\sin2\phi)},\\
B&=&\frac{2(\lambda\cos\phi+x\sin\phi)+A}{\alpha+1},
\eeas
and the first term approximation can be written as
$$
\pol=
L_{n}^{(\alpha)}(\xi)+{{\cal O}}\left(\asp^{n-2}\right),
\quad \xi=A/B,
$$
as $r\to\infty$, uniformly with respect to $\theta$.

As an alternative, we solve $c_1=0, c_2=0$ for 
$A$ and $C$, with  $B=1$. This gives
\beas
A&=&2[x(\sin\phi-\sin2\phi)+\lambda(\cos\phi-\cos2\phi)],\\
C&=&2[x(2\sin\phi-\sin2\phi)+\lambda(2\cos\phi-\cos2\phi)]-1.
\eeas
and the first term approximation can be written as
$$
\pol=
\left[L_{n}^{(\alpha)}(\xi)+{{\cal O}}\left(\asp^{n-2}\right)\right],
\quad \xi=A,\quad \alpha=C.
$$
as $r\to\infty$, uniformly with respect to $\theta$. 

Solving $A=\xi, C=\alpha$ for $x$ and $\lambda$, we obtain
\bea
\lambda&=&(1-\cos\phi)\xi+{1\over2}(\alpha+1)(2\cos\phi-1),\nonumber\\
&{\ }& \\ \label{f8}
x&=&{2(\xi-\alpha-1)\cos^2\phi+
(\alpha+1-2\xi)\cos\phi+\alpha+1-\xi\over 2\sin\phi}.\nonumber
\eea
Then $c_3={2\over3}(\alpha+1-2\xi)(1-\cos\phi)$ and 
$c_k={\cal O}(\phi^2)$ as $\phi\to 0$, 
which follows from deriving a recursion relation
for $c_k$. 

Using these values of $x$ and $\lambda$, we obtain the limit
$$\lim_{\phi \to 0}\pol=L_{n}^{(\alpha)}(\xi).$$

\subsection{Three free parameters}

We solve $c_1=0, c_2=0, c_3=0$ for $A, B$ and $C$. This gives
\beas
A&=&{2\sin\phi(x\sin\phi+\lambda\cos\phi)\WW\over
x\sin2\phi+\lambda\cos2\phi+\sin\phi\WW}         \hskip 4mm
={2r\sin\phi\sin{1\over2}(\theta+\phi)\over\sin{1\over2}(\theta+3\phi)},\\
B&=&{x\sin2\phi+\lambda\cos2\phi+\sin\phi\WW\over
x\sin\phi+\lambda\cos\phi}                       \hskip 4mm
={\sin{1\over2}(\theta+3\phi)\over\sin{1\over2}(\theta+\phi)},\\
C+1&=&2{x\sin2\phi+\lambda\cos2\phi+2\sin\phi\WW\over B^2}  
={2r[\sin(\theta+2\phi)+2\sin\phi]\over B^2}.
\eeas
The first coefficients
$c_k$ are given by
\beas
c_0&=&1,\quad c_1=c_2=c_3=0,\\
c_4&=&{r\over2}\{\sin(\theta+4\phi)+[\sin\phi-\sin(\theta+2\phi)]B^2\}.
\eeas

The first term approximation can be written as
$$
\pol=
B^n\left[L_{n}^{(C)}(\xi)+{{\cal O}}\left(\asp^{n-3}\right)\right],
\quad\xi={A\over B}={2r\sin\phi\over B^2},
$$
as $\asp\to\infty$, uniformly with respect to $\theta$. 

\def\asp{\gamma}
\def\pol{P_n^{(\alpha,\beta)}(x)}
\section{Jacobi, Meixner and Krawtchouk  to Laguerre}\label{sectfive}

We give the results for one free parameter only.

\subsection{Jacobi to Laguerre}

Let $R(w)=\sqrt{1-2xw+w^2}$.
The generating function reads
$$F(x,w)={2^{\alpha+\beta}(1+R-w)^{-\alpha}(1+R+w)^{-\beta}\over R}=
\sum_{n=0}^\infty \pol w^n.
$$
As in Lemma 3.1, we define coefficients $c_k$, and the expansion reads
$$
\pol=
\sum_{k=0}^n B^{n-k}c_k
L_{n-k}^{(C)}(\xi),\quad \xi=A/B.
$$
We consider $\alpha+\beta \to\infty$, and
solve $c_1=0$ for $A$, with  $B=1$ and $C=\alpha$. This gives
$$A=\frac12(\alpha+\beta+2)(1-x).$$
The first coefficients
$c_k$ are given by
$$c_0=1,\quad c_1=0,\quad
c_2=\frac18[-\alpha+3\beta-2(\alpha+3\beta+4)x+(3\alpha+3\beta+8)x^2].
$$ 
The first term approximation can be written as
$$
\pol=
L_{n}^{(C)}(\xi)+{{\cal O}}\left(\asp^{n-1}\right),
\quad  \asp=\alpha+\beta,\quad
\xi=\frac12(\alpha+\beta+2)(1-x).
$$
A limit can be obtained by writing $x=1-2\xi/(\alpha+\beta+2)$ . 
Then we have
$c_k={\cal O}(1/\beta)$ as $\beta\to \infty$ for $k\ge 2$, and we obtain
$$
\lim_{\beta\to\infty}P_n^{(\alpha,\beta)}[1-2\xi/(2+\alpha+\beta)]=
L_{n}^{(\alpha)}(\xi),
$$
which includes the limit of the Askey scheme (cf. \cite{koeksw})
$$
\lim_{\beta\to\infty}P_n^{(\alpha,\beta)}(1-2\xi/\beta)=
L_{n}^{(\alpha)}(\xi).
$$
\def\asp{\beta}
\def\pol{M_n(x;\beta,c)}

\subsection{Meixner to Laguerre}

The generating function reads
$$F(w)=\left(1-{w\over c}\right)^{x}
\left(1-w\right)^{-\beta-x}=
\sum_{n=0}^\infty {(\beta)_n\over n!}M_n(x;\beta,c)w^n,
$$
and we define $c_k$ as in Lemma 3.1.
The expansion reads
$$
\pol=
\sum_{k=0}^n B^{n-k}c_k
L_{n-k}^{(C)}(\xi),\quad \xi=A/B.
$$
We solve $c_1=0$ for $A$, with  $B=1$ and $C=\alpha$. This gives
$$A={(\alpha-\beta+1)c+(1-c)x\over c}.$$
The first coefficients
$c_k$ are given by
$$
c_0=1,\quad c_1=0,\quad
c_2={(1+\alpha-\beta)c^2+(2c-c^2-1)x\over 2c^2}.
$$
The first term approximation can be written as
$$
\pol=
L_{n}^{(\alpha)}(\xi)+{{\cal O}}\left(\asp^{n-1}\right),\quad 
\xi={(\alpha-\beta+1)c+(1-c)x\over c}.
$$
A limit can be obtained by putting $\beta=\alpha+1$ 
and writing $x=c\xi/(1-c)$. 
Then we have $c_2=(c-1)\xi/(2c)$, and
$c_k={\cal O}(1-c)$ as $c\to 1$ for $k\ge 2$. We obtain
the limit of the Askey scheme (cf. \cite{koeksw})
$$
\lim_{c\to 1}M_n(c\xi/(1-c);\alpha+1,c)=
{L_{n}^{(\alpha)}(\xi)\over L_{n}^{(\alpha)}(0)}.
$$

\def\asp{N}
\def\pol{{N\choose n}K_n(x;p,N)}

\subsection{Krawtchouk to Laguerre}

Let 
$q:={1-p\over p}$. The generating function reads
$$F(w)=\left(1-qw\right)^{x}
\left(1+w\right)^{N-x}=
\sum_{n=0}^N \pol w^n,
$$
and we define $c_k$ as in Lemma 3.1. The expansion reads
$$
\pol=
\sum_{k=0}^n B^{n-k}c_k
L_{n-k}^{(C)}(\xi),\quad \xi=A/B.
$$
We solve $c_1=0$ for $A$, with  $B=1$ and $C=\alpha$. This gives
$A= \alpha + 1-N + (1+q)x$.
The first coefficients
$c_k$ are given by
$$c_0=1,\quad c_1=0,\quad c_2={1\over2}[1+\alpha-3N+(3+2q-q^2)x].$$
The first term approximation can be written as
$$
\pol=
L_{n}^{(C)}(\xi)+{{\cal O}}\left(\asp^{n-1}\right),\quad 
\xi=\alpha + 1-N + (1+q)x,\quad \asp \to \infty.
$$

\end{document}